\documentclass[12pt]{article}
\usepackage{amsmath}
\usepackage{graphicx}
\usepackage{geometry}
\usepackage{lipsum} 
\usepackage{indentfirst}

\usepackage{hyperref}
\usepackage{xcolor}
\usepackage{graphicx}
\usepackage{subcaption}
\usepackage{float}
\usepackage[utf8]{inputenc}
\usepackage{tikz}
\definecolor{krpurple}{RGB}{128,0,128} 

\definecolor{Aqua}{rgb}{0.2, 1.0, 1.0} 
\title{Inverse Design of a Layered Medium for Optimal Surface Localization\thanks{This is a preprint of the article published in SIAM Undergraduate Research Online (SIURO) Volume 18, 2025. The published version is available at \url{https://doi.org/10.1137/24S1713685}.}}
\author{Ziling Chen $^*$\\ Project Advisor: Fadil Santosa $^\dagger$}
\date{}

\begin{document}

\maketitle

\begin{abstract}
    \noindent Electromagnetic wave manipulation plays a crucial role in advancing technology across various domains, including photonic device design. This study presents an inverse design approach for a periodic medium that optimizes electromagnetic wave localization at the interface between a layered half-space and a homogeneous half-space. The approach finds a maximally localized mode at a specified frequency and wave number. The mode propagates in the direction of the interface. The design parameters are the permittivity of the layered medium, their relative thicknesses, and the permittivity of the homogeneous half-space. We analyze the problem using the transfer matrix method and apply the particle swarm optimization to find a rapidly decaying mode that satisfies the design constraints. The design process is demonstrated in a numerical example, which serves to illustrate the efficacy of the proposed method. 
\end{abstract}

\vfill 

\begin{flushleft}
\textbf{Author Note} \\
Ziling Chen, Department of Applied Math and Statistics, Johns Hopkins University, Email: \href{mailto:zchen105@jh.edu}{zchen105@jhu.edu}. \\
Fadil Santosa, Advisor, Department of Applied Math and Statistics, Johns Hopkins University, Email: \href{mailto:fsantos9@jh.edu}{fsantos9@jhu.edu}.
\end{flushleft}

\newpage
\section{Introduction}

Wave localization\cite{Figotin1997} confines waves -- light, sound, or other forms -- to a small region within a medium, countering their natural tendency to spread. This phenomenon is often achieved through specific structural designs. A photonic band gap\cite{PhysRevLett.58.2059}\cite{yablonovitch1993} is a frequency range where electromagnetic waves cannot propagate in any direction within a medium. Within this band gap, optical modes, spontaneous emission, and zero-point fluctuations are absent.

Previous study by Figotin and Gorentsveig\cite{PhysRevB.58.180} extensively explored the method aspects of wave localization through two-layered media and the effect of the defects on the localization, providing an illustration of the transfer matrix method and localized mode pattern at varying medium configurations. 
A more recent study by Rammohan\cite{Rammohan2010} explores the application of a hybrid optimization algorithm that combines the bio-inspired genetic algorithm\cite{Holland1992} with gradient descent. This approach aims to achieve global optimum search and precise manipulation of photonic band gaps.

In our study, we apply the transfer matrix method for the band gap calculation\cite{PhysRevLett.69.2772}\cite{pendry1996TMM}\cite{pendry1996Bandgap} to model wave localization in a layered half-space joined to a homogeneous half-space. The layered structure is periodic with each period consisting of two homogeneous layers. The model parameters are the relative thicknesses of the two constituents and the material constant of the homogeneous half-space.  Our goal is to solve the following design problem: For specified frequency and wave number (along the direction of propagation), find the model parameters that produces a mode that is maximally localized at the interface between the two half-spaces.  To solve the optimization problem, we use a bio-inspired optimization algorithm, particle swarm optimization (PSO)\cite{Kennedy1995}.

In Section 2 of this paper, we introduce the assumptions and medium under consideration, followed by the derivation of the wave equation, which serves as the mathematical model. Section 3 then elaborates on the transfer matrix method derived from the governing equation. Section 4 discusses the selection of parameters for the localized mode in a homogeneous half-space and the conditions necessary for achieving localization\cite{PhysRevB.58.180}\cite{Schultz1993}\cite{PhysRevA.39.2005}. The main ideas are illustrated in a numerical example described in the same section. Section 5 employs particle swarm optimization (PSO) to address the inverse problem, optimizing the design parameters. The results from this optimization are presented in the same section. In section 6, we provide a summary of our findings and indicate possible future research directions.

\section{Model of the problem}
\begin{figure}[htp]
    \centering
    \begin{tikzpicture}[scale = 0.7]
    \draw[->] (-2,0) -- (11,0) node[right] {$x$};
    \draw[->] (0,0) -- (0,6) node[above] {$y$};

    \fill[yellow!10] (-2,0) rectangle (1,5);
    \fill[blue!20] (0,0) rectangle (1,5);
    \fill[blue!20] (3,0) rectangle (4,5);
    \fill[blue!20] (6,0) rectangle (7,5);
    \fill[blue!20] (9,0) rectangle (10,5);
    \node[rotate=90, scale=0.7] at (-1,2.5) {Homogeneous Half Space};
    \node at (0.5,3) {A};
    \node[below] at (0,0) {$0$};
    \node at (2,3) {B};
    \node[below] at (3,0) {$p$};
    \node at (3.5,3) {A};
    \node at (5,3) {B};
    \node[below] at (6,0) {$2p$};
    \node at (6.5,3) {A};
    \node at (8,3) {B};
    \node[below] at (9,0) {$3p$};
    \node at (9.5,3) {A};

    \draw[<->] (0,0.5) -- (3,0.5) node[midway, above] {$p$};
    \draw[<->] (3,0.5) -- (6,0.5) node[midway, above] {$p$};
    \draw[<->] (6,0.5) -- (9,0.5) node[midway, above] {$p$};
\end{tikzpicture}
\captionsetup{font = small}
    \caption{For $x>0$, the medium is a periodically layered medium where each period consists of layers of materials A and B, the period is $p$. For $x<0$, the medium is homogeneous.
 }
    \label{fig:medium_visualization}
\end{figure}
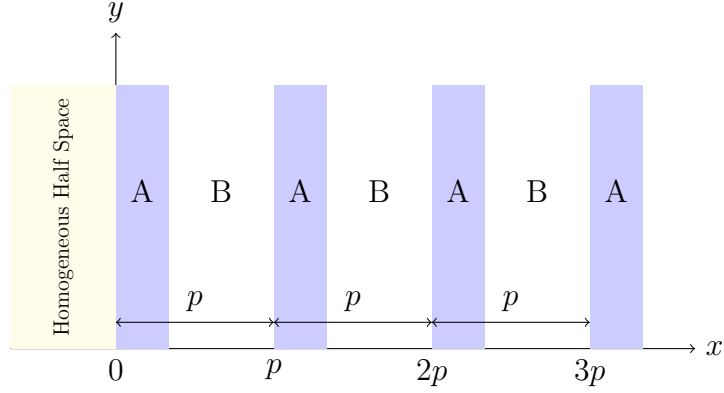
In this work, we consider the propagation of electromagnetic waves in a two-dimensional heterogeneous medium. Using coordinates $xy$, we place a homogeneous half-space in $x<0$. For $x>0$, the medium is layered in the $x$-direction with periodic units of width $p$ (see Figure \ref{fig:medium_visualization}).  We will assume a Transverse Magnetic (TM) field, meaning that the electric vector field $\mathbf{E}$ and the magnetic vector field $\mathbf{H}$ are of the form
\[
\mathbf{E} = (E_x, E_y, 0), \quad \mathbf{H} = (0, 0, H_z).
\]
All dependent variables are functions of $(x,y)$ and time $t$. Maxwell's equations governs the behavior of $\mathbf{E}$ and $\mathbf{H}$
\begin{align}
\mu \frac{\partial \mathbf{H} }{\partial t} &= \nabla \times \mathbf{E}, \label{Maxwell_1}\\
 \epsilon \frac{\partial \mathbf{E} }{\partial t} &= \nabla \times \mathbf{H}, 
\label{Maxwell_2}
\end{align}
where $\epsilon$ and $\mu$ are the electric permittivity and magnetic permeability of the medium. In our model, we will assume that $\mu$ is constant and $\epsilon$ is a function of $x$. Using the form of TM modes and their functional dependence, Maxwell's equation \eqref{Maxwell_1}-\eqref{Maxwell_2} reduces to
\begin{equation}
\begin{aligned}
\mu \frac{\partial H_z}{\partial t} &= \frac{\partial E_x}{\partial y}-\frac{\partial E_y}{\partial x}, \\
\epsilon \frac{\partial E_x}{\partial t} &= \frac{\partial H_z}{\partial y}, \\
\epsilon \frac{\partial E_y}{\partial t} &= -\frac{\partial H_z}{\partial x} .
\end{aligned}
\label{derivation_eqs}
\end{equation}
We take the time derivatives of the first equation in \eqref{derivation_eqs} to get
\[
\mu \frac{\partial^2 H_z}{\partial t^2} = -\frac{\partial^2 E_y}{\partial t \partial x} + \frac{\partial^2 E_x}{\partial t \partial y}.
\]
Next, we take time derivatives of the second and third equations in \eqref{derivation_eqs} and substitute them in the above to obtain a single equation for $H_z$
\begin{equation}
\epsilon \mu \frac{\partial^2 H_z}{\partial t^2} = \frac{\partial^2 H_z}{\partial x^2} + \frac{\partial^2 H_z}{\partial y^2}.
\label{Hz-equation}
\end{equation}
Since the material properties are only $x$-dependent, we can assume a separable solution of the form
\[
H_z(x, y, t) = u(x) e^{i(\eta y - \omega t)},
\]
where \(\eta\) is the spatial ($y$) wave number, and \(\omega\) is the angular time frequency. Once this form is substituted in \eqref{Hz-equation} we obtain a differential equation for $u(x)$. 
\[
u''(x) + (\mu \epsilon(x) \omega^2 - \eta^2) u(x) = 0 .
\]
We will use the speed of light
\[
c(x) = \frac{1}{\sqrt{\mu \epsilon(x)}},
\]
as the function describing the medium.  The final governing equation for $u(x)$ is
\begin{equation}
    u''(x) + \left(\frac{\omega^2}{c(x)^2}-\eta^2 \right) u(x) = 0.
    \label{WaveEQ}
\end{equation}
The remainder of this work will focus on finding closed-form solutions for $u(x)$.
 
The function $c(x)$ can now be made precise.  For $x<0$, $c(x) = c_0$ where $c_0$ is a constant that we will determine from our calculation. For $x>0$, $c(x)$ is a periodic function where within each period $p$ we have
\begin{equation}
    c(x) = \left\{ \begin{array}{cc}
                c_A , & 0 \leq x < \theta p \\
                c_B , & \theta p \leq x < p 
                   \end{array} \right. ,
\end{equation}
where $\theta$ represents the volume fraction of material $A$ in each periodic layer.
\section{Solution for $u(x)$}
Since the medium is piecewise constant in $x$, we define $\sigma_0 = \sqrt{\frac{\omega^2}{c_0^2}-\eta^2}$. Then, for $x < 0$, 
\begin{equation}
    u(x) = C_1e^{i\sigma_0x}+C_2e^{-i\sigma_0x},
    \label{u_x_neg}
\end{equation}
For $x > 0$, we use a trigonometric form to solve \eqref{WaveEQ}
\begin{equation}
    u(x) = A\cos(\sigma x) + B\sin(\sigma x),
    \label{ux_1}
\end{equation}
\begin{equation}
    u'(x) = -\sigma A \sin(\sigma x) + \sigma B \cos(\sigma x),
    \label{ux_p}
\end{equation}
where $\sigma = \sigma_A$ or $\sigma_B$ depending on $x$, and $\sigma_A = \sqrt{\frac{\omega^2}{c_A^2}-\eta^2}$, $\sigma_B = \sqrt{\frac{\omega^2}{c_B^2}-\eta^2}$.\\
We normalize the coordinates $x$ such that $p=1$ for simplicity of the computation.

Consider the first periodic cell $0\leq x\leq 1$ (see figure \ref{fig:position_normalization}). For $x=0$, \eqref{ux_1}-\eqref{ux_p} provides the initial condition $\begin{bmatrix}
    u\\
    u'
\end{bmatrix}(0) = \begin{bmatrix}
    A\\
    \sigma_AB
\end{bmatrix}. $ Therefore, for $0<x<\theta$, we have
\begin{equation}
\begin{aligned}
    \begin{bmatrix}
        u \\
        u'
    \end{bmatrix}(x) = 
    \begin{bmatrix}
        \cos(\sigma_A x) & \sin(\sigma_A x) / \sigma_A \\
        -\sigma_A \sin(\sigma_A x) & \cos(\sigma_A x)
    \end{bmatrix}
    \begin{bmatrix}
        u \\
        u'
    \end{bmatrix}(0).
\end{aligned}
\label{wave-propag-eq}
\end{equation}
Writing 
\begin{equation}
    T_A = \begin{bmatrix}
            \cos(\sigma_A \theta) & \sin(\sigma_A \theta) / \sigma_A \\
            -\sigma_A \sin(\sigma_A \theta) & \cos(\sigma_A \theta)
            \end{bmatrix},
\end{equation}
we have
\begin{equation}
    \begin{bmatrix}
        u \\
        u'
    \end{bmatrix}(\theta) = T_A \begin{bmatrix}
        u \\
        u'
    \end{bmatrix}(0).
\end{equation}
For $\theta < x < 1$, we have 
\begin{equation}
\begin{aligned}
    \begin{bmatrix}
        u \\
        u'
    \end{bmatrix}(x) = 
    \begin{bmatrix}
        \cos(\sigma_B (x-\theta)) & \sin(\sigma_B (x-\theta)) / \sigma_B \\
        -\sigma_B \sin(\sigma_B (x-\theta)) & \cos(\sigma_B (x-\theta))
    \end{bmatrix}
    \begin{bmatrix}
        u \\
        u'
    \end{bmatrix}(\theta).
\end{aligned}
\label{wave-propag-eq-B}
\end{equation}
Letting 
\begin{equation}
    T_B = 
            \begin{bmatrix}
            \cos(\sigma_B (1-\theta)) & \sin(\sigma_B (1-\theta)) / \sigma_B \\
            -\sigma_B \sin(\sigma_B (1-\theta))) & \cos(\sigma_B (1-\theta))
            \end{bmatrix},
\end{equation}
we have
\begin{equation}
    \begin{bmatrix}
        u \\
        u'
    \end{bmatrix}(1) = T_B \begin{bmatrix}
        u \\
        u'
    \end{bmatrix}(\theta).
\end{equation}
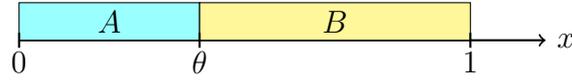
\begin{figure}[H]
    \centering
    \begin{tikzpicture}
    \def\totalLength{6} 
    \def\thetaFraction{0.4} 

    \pgfmathsetmacro{\thetaP}{\totalLength*\thetaFraction}
    \pgfmathsetmacro{\oneMinusThetaP}{\totalLength - \thetaP}
    \draw[thick, ->] (0,0) -- (\totalLength + 1,0) node[anchor=west] {$x$};
    \filldraw[fill=Aqua!50] (0,0) rectangle (\thetaP,0.5);
    \filldraw[fill=yellow!50] (\thetaP,0) rectangle (\totalLength,0.5);

    \node at (\thetaP/2,0.25) {$A$}; 
    \node at (\thetaP + \oneMinusThetaP/2,0.25) {$B$}; 
    \draw[thick] (0,-0.1) -- (0,0.1) node[anchor=north, yshift=-1mm] {0};
    \draw[thick] (\thetaP,-0.1) -- (\thetaP,0.1) node[anchor=north, yshift=-1mm] {$\theta$};
    \draw[thick] (\totalLength,-0.1) -- (\totalLength,0.1) node[anchor=north, yshift=-1mm] {1};
\end{tikzpicture}
\caption{This figure illustrates the combination of the first period $p = 1$ of the medium. The relative thickness of material A is $\theta$ and the relative thickness of material B is $1-\theta$.}
\label{fig:position_normalization}
\end{figure}
\noindent Then we have the wave propagation after the first periodic cell
\begin{equation}
    \begin{bmatrix}
        u \\
        u'
    \end{bmatrix}(1) = T_BT_A \begin{bmatrix}
        u \\
        u'
    \end{bmatrix}(0).
\end{equation}
For our target medium, the transfer matrix \cite{pendry1996TMM} after \(n\) periods (\(x = n\)) is defined as
\begin{equation}
T^{(1)} = (T_B T_A)^n,
\end{equation}
satisfying:
\begin{equation}
\label{TMM}
    \begin{bmatrix}
        u \\
        u'
    \end{bmatrix}(n) = T^{(1)} \begin{bmatrix}
        u \\
        u'
    \end{bmatrix}(0).
\end{equation}
\section{Localized mode}
\subsection{Decay in Medium}
For our problem, a localized mode is characterized by an exponential decrease in the amplitude as $x$ increases. That is, it refers to a solution $u(x)$ that decays away from $x=0$ in the positive $x$ direction. Such a mode can be found by examining the eigenvalues of $T^{(1)}$. Let us write 
\begin{equation}
    [T_BT_A] = V\Lambda V^{-1},
\end{equation}
then 
\begin{equation}
    T^{(1)} = V\Lambda^{n} V^{-1}.
\end{equation}
By the dimensions of $T_A$ and $T_B$, $T^{(1)}$ is a $2\times 2$ matrix and has two eigenvalues, denoted as $\lambda_1$ and $\lambda_2$ with their associated eigenvectors $v_1$ and $v_2$. 
Since all entries of $T^{(1)}$ are real, both eigenvalues of $T^{(1)}$ are either complex conjugate or real. Note that since $T_A$ and $T_B$ are unitary matrices, the eigenvalues are mutually inverse \cite{PhysRevB.58.180}, that is, $\lambda_1\lambda_2 = 1$. A localized mode is one for which $\lambda_1$ and $\lambda_2$ are both real and $|\lambda_1| < 1$. Moreover, the rate of decay in the amplitude of $u(x)$ is inversely proportional to $|\lambda_1|$. 

\subsection{Decay in Negative Homogeneous Half-Space}\label{DecayInNegHomSpace}
In the case of negative $x$, for the solution $u(x)$, to decay, we choose $c_0$ such that the coefficient before $x$ in \eqref{u_x_neg} is positive and $C_1$ in the equation is 0. This means that
\begin{equation}
    \frac{\omega^2}{c_0^2}-\eta^2 < 0,
    \label{neg-dir-decay-cond}
\end{equation}
and the solution $u(x)$ becomes
\begin{equation}
    u(x) = C_2e^{-i\sigma_0 x}.
    \label{ux_neg}
\end{equation}
The solution vector at $x=0$ can now be derived, and is given by
\begin{equation}
    \begin{bmatrix}
        u\\
        u'
    \end{bmatrix}(0) = C_2\begin{bmatrix}
        1\\
        -i\sigma_0
    \end{bmatrix}
\end{equation}
To determine the value of $c_0$, we align the eigenvector associated with the smaller absolute eigenvalue, \textbf{$v_1$} with the initial solution vector $\begin{bmatrix}
        u \\
        u'
    \end{bmatrix}(0)$ to ensure the continuity of values of $u(x)$ and $u'(x)$ at the point $x=0$. Mathematically, this means that
    \begin{equation*}
        C_2\begin{bmatrix}
        1\\
        -i\sqrt{\frac{\omega^2}{c_0^2}-\eta^2}
    \end{bmatrix} = \alpha v_1, \;\; \mbox{for} \;\; \alpha \in \Re.
    \end{equation*}
Let $v_{11}$ and $v_{21}$ denote the first and second entry of $v_1$ respectively. Define $r = v_{11}/v_{21}$ and setting $\alpha = C_2$, we have
    \begin{equation*}
        \frac{1}{-i\sqrt{\frac{\omega^2}{c_0^2}-\eta^2}} = r .
    \end{equation*}
Therefore, we must choose
    \begin{equation}
        c_0^2 = \frac{\omega^2}{\eta^2-\frac{1}{r^2}}.
        \label{c0cond}
    \end{equation}
Here $c_0$ must be real and positive. We can see that $c_0$ depends on the eigenvectors of the transfer matrix. 

\subsection{An Example}
\begin{figure}[H]
    \centering
    \begin{subfigure}[b]{0.4\textwidth}
        \centering
        \includegraphics[width=\linewidth]{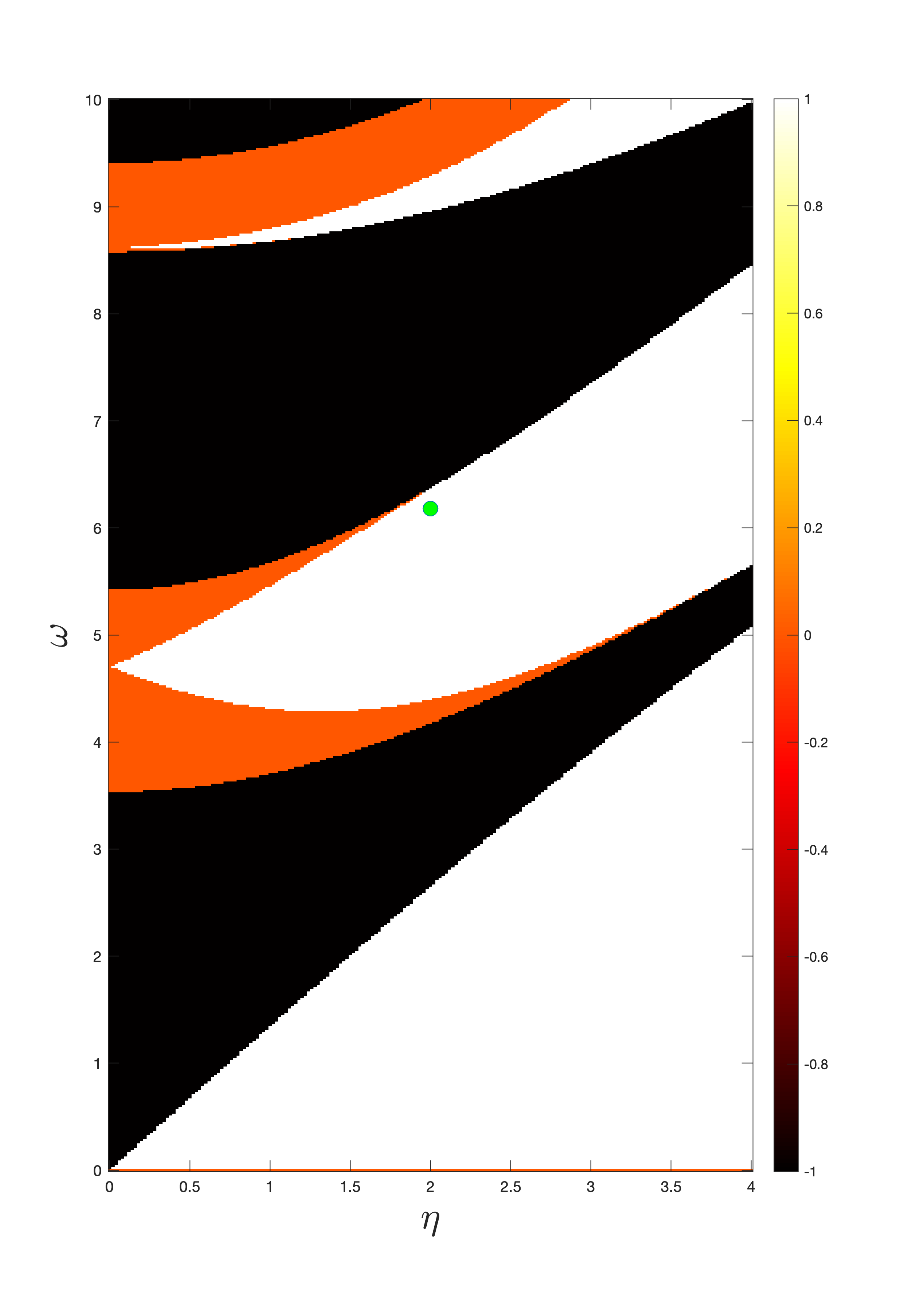}
        \caption{}
        \label{fig:band_gaps}
    \end{subfigure}
    \begin{subfigure}[b]{0.55\textwidth}
        \centering
        \includegraphics[width=\linewidth]{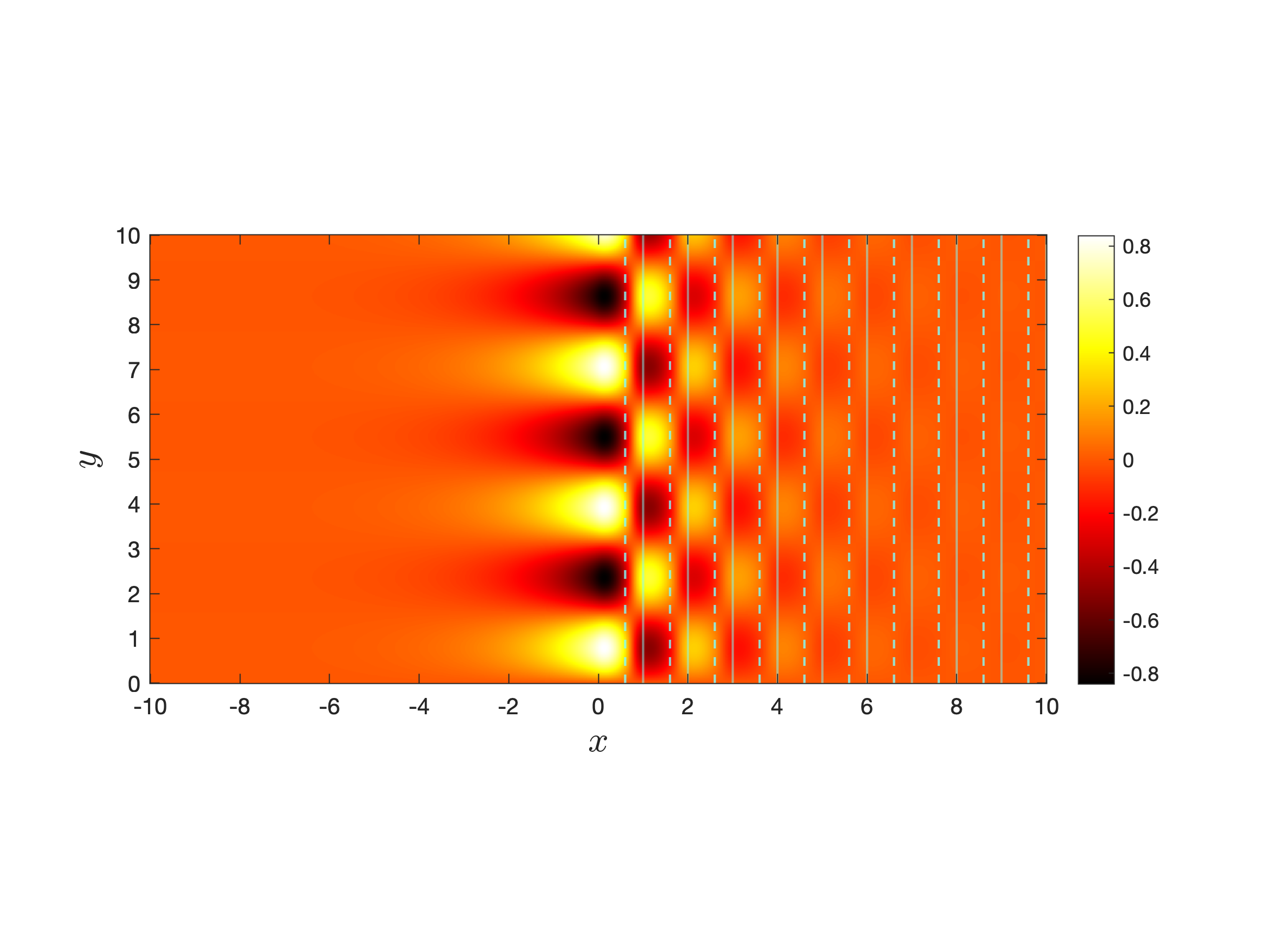}
        \caption{}
        \label{fig:wave_propagation_complete}
    \end{subfigure}
    \captionsetup{font=small}
   \caption{We fix the medium parameters $c_A = 2$, $c_B = 1$, and $\theta = 0.6$. (a) This subfigure illustrates the wave localization behavior across values of $(\eta, \omega)$. The orange and white regions indicate regions in which waves decay in $x>0$, with the white specifically denoting areas where the waves also decay in the negative direction, thereby localizing at $x=0$. A green circle highlights the specific parameters $\omega = 6.18$ and $\eta = 2$, under which the conditions for wave localization are met for the chosen medium combination. (b) This subfigure visualizes the wave amplitude $u(x)\sin \eta y$ corresponding to the parameters marked by the green circle in subfigure (a). In the right proportion of this subfigure where $x>0$, we use solid and dashed lines to exhibit the medium structure. Solid lines mark the boundaries of periodic cells, while dashed lines differentiate materials, indicating their respective volume fractions ($\theta = 0.6$ in this example). This subfigure effectively shows the wave localization at $x=0$.}
    \label{fig:band_gap_wave_propag}
\end{figure}

We can now examine the direct problem. We start by setting $c_A$, $c_B$, and $\theta$. In this section, we fix the medium parameters $c_A = 2$, $c_B = 1$, and $\theta = 0.6$. Next, we compute the eigenvalues of the matrix $[T_B T_A]$ for a range of values of $(\eta,\omega)$. Figure \ref{fig:band_gaps} illustrates the choice of wave parameters $(\eta,\omega)$ that satisfy the wave localization. The white region with a value of 1 in this figure means the wave decays in both $x$ directions. This means that for $(\eta,\omega)$ in the white region, $|\lambda_1|<1$ and $c_0$ in \eqref{c0cond} is real. In the orange region, we have $|\lambda_1|<1$, but we cannot find a real $c_0$ from \eqref{c0cond}. In the black region, the waves fail to decay in the $x>0$.

After choosing a pair of wave parameters, $\omega = 6.18$ and $\eta = 2$ (green circle in figure \ref{fig:band_gaps}) satisfying the conditions mentioned in previous parts in this section, with the value of $|\lambda_1|= 0.5996 < 1$.  We then compute $u(x)\sin \eta y$ for the wave propagation visualization. Figure \ref{fig:wave_propagation_complete} gives an example of the wave propagating in both the negative homogeneous half-space and our medium after 10 periods. It can be seen that the wave propagates in the $y$-direction and decays in both positive and negative $x$ directions.
\section{Inverse Design Problem}
The study in the previous section clearly shows that the rate of decay of the wave for $x>0$ depends on the properties of the layered medium. We now address the inverse design problem of determining medium properties that maximally localize the wave. 
\subsection{Objective Function and Variables}
The primary variables in this inverse problem are the speed of light in the materials, \(c_A\) and \(c_B\), and the volume fraction \(\theta\). Let $\lambda_1$ be the eigenvalue of $[T_BT_A]$ such that $|\lambda_1|\leq1$. The decay rate into $x>0$ is associated with $\lambda_1$. To make the decay rapid, we must make $\lambda_1$ to be small. To this end, we define 
\begin{equation}
f(c_A, c_B, \theta) = 
    \begin{cases}
    2, &|\lambda_1|=1\\
    |\lambda_1|, &|\lambda_1|<1
    \end{cases}.
    \label{objectiv_fct}
\end{equation}
Minimizing $f(c_A, c_B, \theta)$ provides a medium configuration for which a wave with wave number $\eta$ and frequency $\omega$ decays the fastest in the medium.
\subsection{Particle Swarm Optimization}
Given the discontinuity and potential for multiple local minima in the landscape of the objective function, an optimization model capable of exploring a broad set of initial conditions is preferred. It has been demonstrated by \cite{Wetter2004} that particle swarm optimization (PSO) is robust and well-suited for discontinuous objective functions. PSO simulates a population of potential solutions that collectively explore the search space, adjusting their positions based on both their locally best-found positions and the global optimal among them. This search mechanism ensures that the search for an optimal material combination is exhaustive and less likely to be trapped in suboptimal solutions, thereby increasing the likelihood of finding a global optimum that meets the desired wave propagation criteria.
\subsection{Optimization Model Formulation}
For simplicity of comparison, we fix the target wave frequency and wave number established in Section 4 ($\omega = 6.18$ and $\eta = 2$). Several constraints are imposed to maintain the physical correctness and feasibility of the solution. Firstly, two different effective wave speeds \(c_A\) and \(c_B\) are bounded within the range \([0.5, 3.5]\). Then, for each $c_0$ calculated from each transfer matrix of the medium, the condition \eqref{neg-dir-decay-cond} has to be satisfied. Additionally, in each iteration, we calculate $r$ as stated in Section \ref{DecayInNegHomSpace}, and we must enforce that \eqref{c0cond} produces $c_0$ that is positive.

To ensure an efficient decay of electromagnetic waves in $x<0$, we add a condition on the exponential coefficient in Equation \ref{ux_neg} that the coefficient should be large enough to ensure the amplitude of the wave decays to 0 after 10 periods of propagation. 

With all the conditions, the optimization model of the inverse problem can be formulated as follows:
\begin{align*}
\text{min} \quad & f(c_A, c_B, \theta) \\
\text{s.t.: } &
 c_A \neq c_B \\
 &\sigma_A>0\\
 &\sigma_B>0\\
&\frac{\omega^2}{c_0^2}-\eta^2 < 0\\
&-i\sigma_0 >0.5\\
&c_0^2 = \frac{\omega^2}{\eta^2-\frac{1}{r^2}}>0\\
&c_A \in [0.5, 3.5], c_B \in [0.5, 3.5], \theta \in (0, 1)
\end{align*}

\subsection{Results}

\begin{figure}[H]
    \centering
    \includegraphics[width=0.60\linewidth]{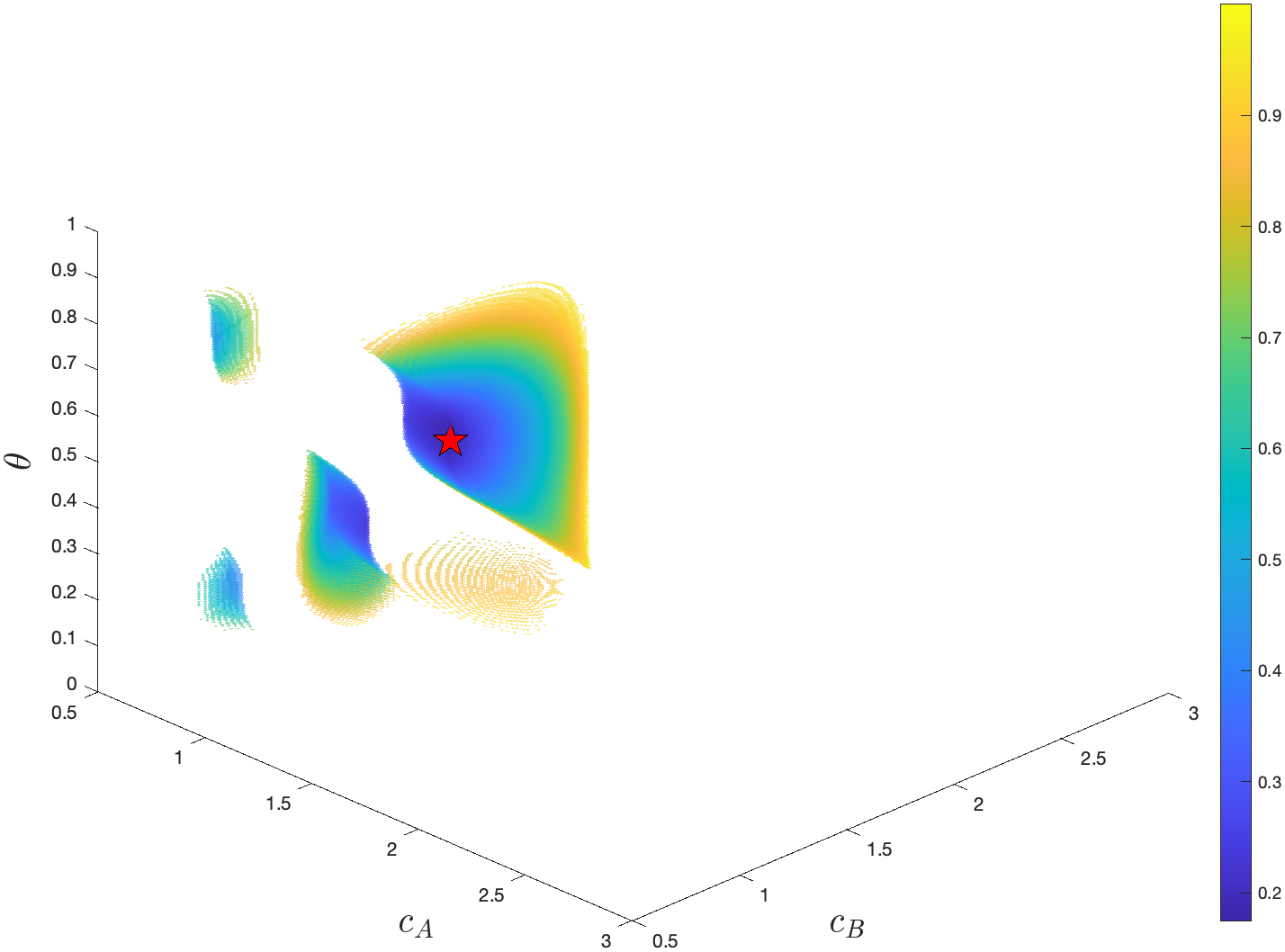}
    \captionsetup{font = small}
    \caption{This figure illustrates the objective function values within a three-dimensional solution space, where the optimal solution is highlighted by a red star. The red star is located in the deepest color region, visually affirming its status as the most optimal point across the entire solution space. This 3D representation effectively underscores the thoroughness of the optimization analysis.}
    \label{fig:optimality_3d}
\end{figure}
The optimal solution given by PSO is \(c_A = 2.15, c_B = 0.50\), and \(\theta = 0.87\), found after 90 iterations. The objective function value is now $|\lambda_1| = 0.1742$. Figure \ref{fig:optimality_3d} presents a 3D graph visualizing the objective function value in the feasible solution space. The blank region indicates that the medium combination corresponding to the solution point fails to cause wave decays in both $x$-directions. The colored region in the figure represents the objective function values of the feasible solution, and the red star highlights the optimal solution found by the PSO. 

\begin{figure}[H]
    \centering
    \begin{subfigure}[b]{0.4\textwidth}
        \centering
        \includegraphics[width=\linewidth]{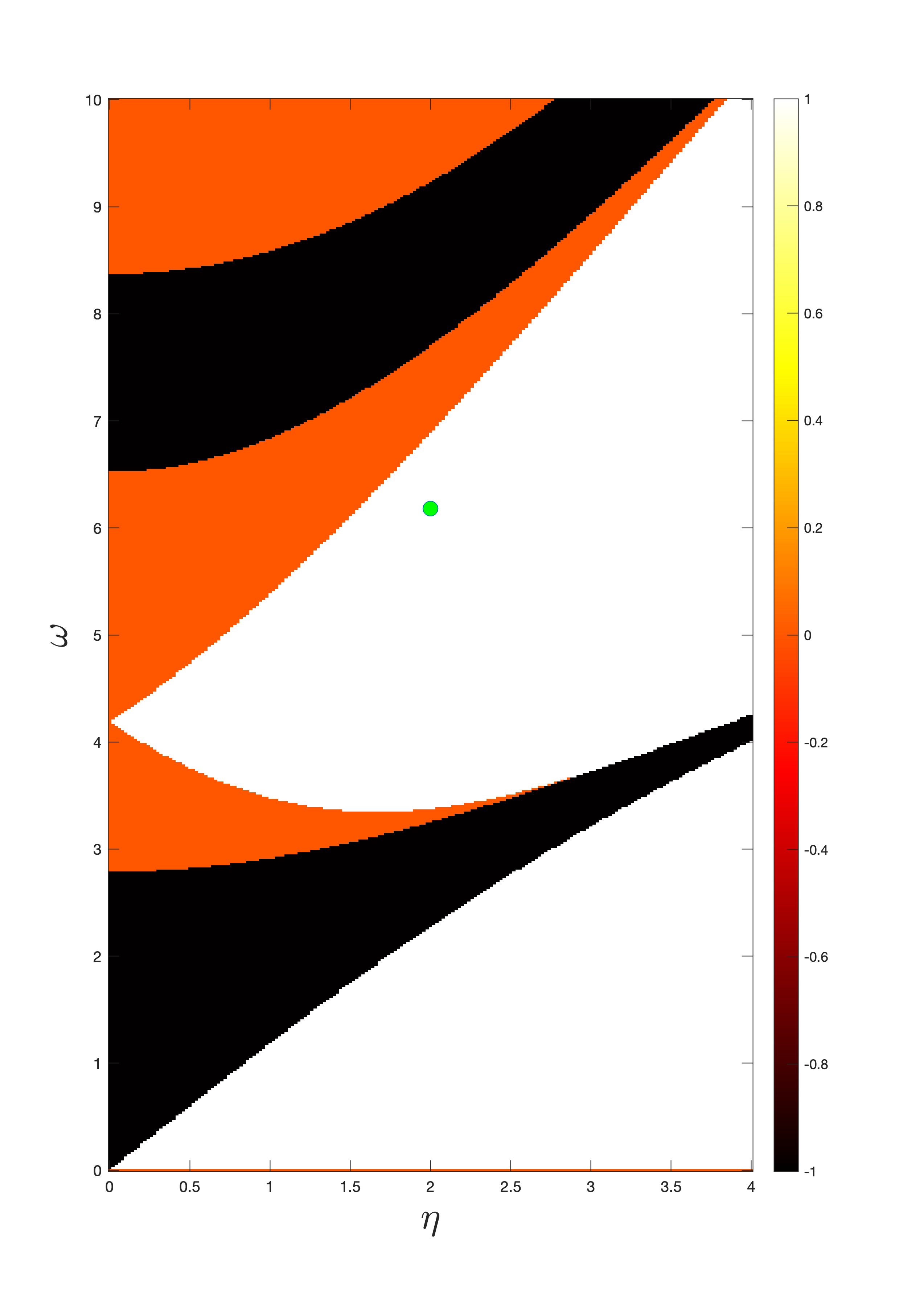}
        \caption{}
        \label{fig:band_gaps_opt}
    \end{subfigure}
    \begin{subfigure}[b]{0.55\textwidth}
        \centering
        \includegraphics[width=\linewidth]{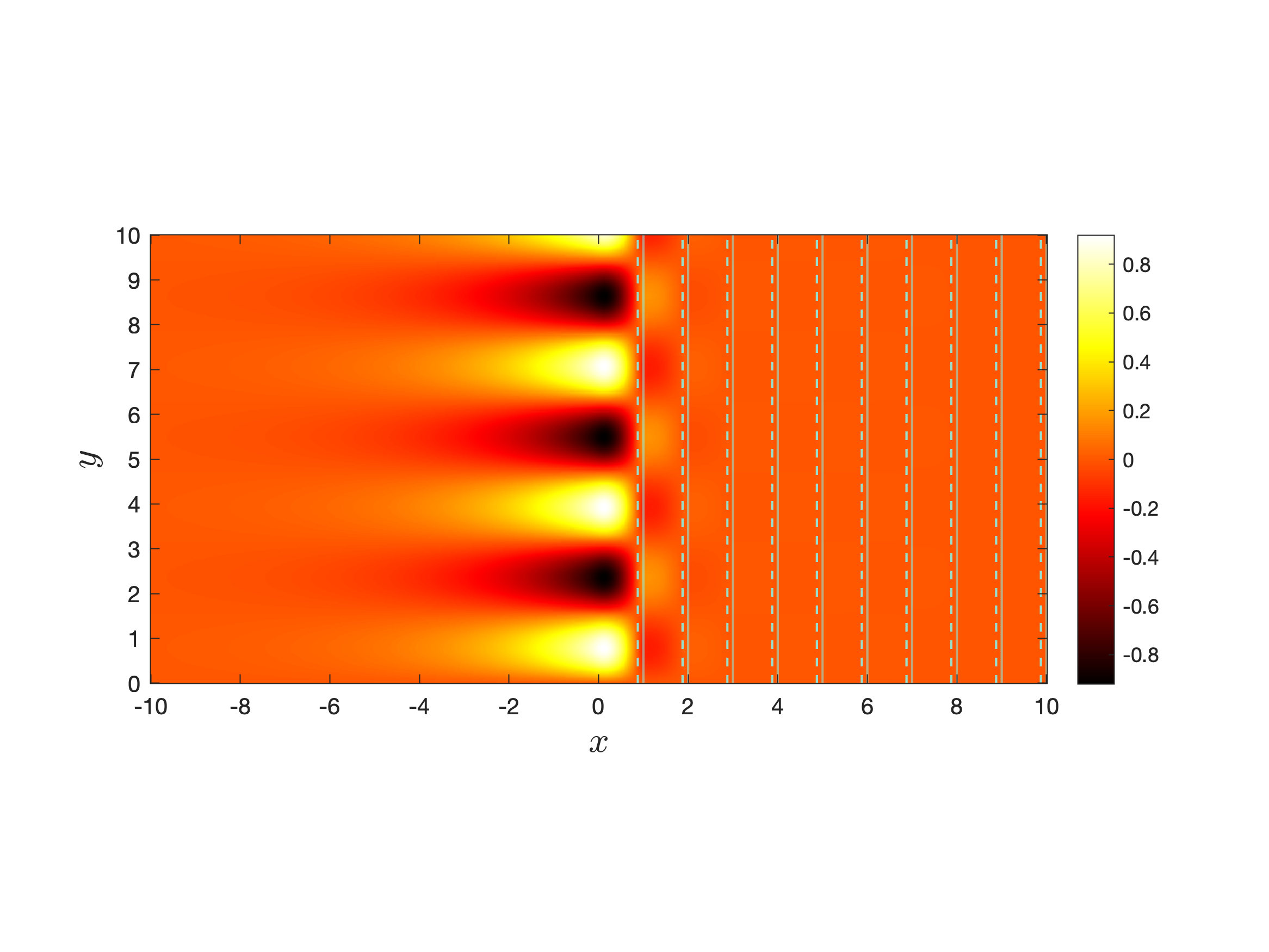}
        \caption{}
        \label{fig:wave_propagation_opt}
    \end{subfigure}
    \caption{(a) This subfigure displays the chosen wave parameters $\omega$ and $\eta$ for the optimized medium composition of $c_A = 2.15$, $c_B = 0.50$, and $\theta = 0.87$, denoted by a green circle in the white zone, indicating wave localization at $x=0$. (b) This figure visualizes wave propagation within this optimized medium configuration, facilitating an easy comparison with previous images and showing accelerated decay. The solid lines and dashed lines also exhibit the structure of the optimized medium, with dashed lines marking materials with optimized volume fraction $\theta = 0.87$.}
    \label{fig:wave_propag_opt}
\end{figure}

\begin{figure}[H]
    \centering
    \includegraphics[width=0.6\linewidth]{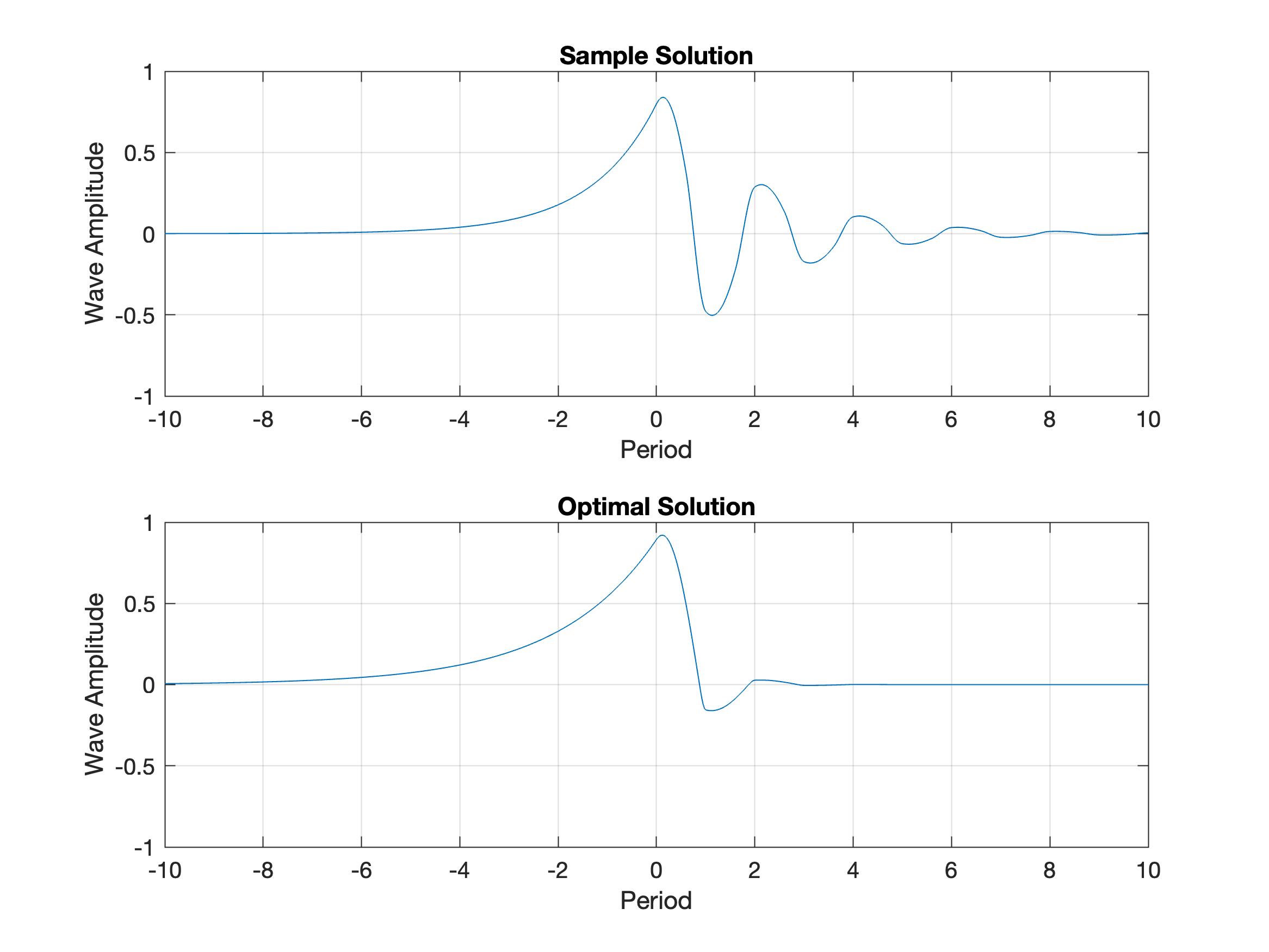}
    \caption{This figure provides a direct comparison of wave amplitude changes over periods with different material combinations. The wave amplitude decays to zero after approximately the $3_{rd}$ period, while in the original medium, the wave decays to near zero after the $8_{th}$ period.}
    \label{fig:wave_amplitude_comparison}
\end{figure}

The wave propagation in the designed medium is visualized in Figure \ref{fig:wave_propag_opt}. 
Figure \ref{fig:wave_amplitude_comparison} provides a comparison between the wave amplitude changes over periods. It can be seen from the figure that the wave in our designed medium decays more rapidly in both directions and localizes at $x = 0$.

\section{Discussion}
The purpose of this study is to demonstrate that it is possible to design a layered medium for maximal wave localization. We use the transfer matrix method to characterize the solution to the wave equation and pose a constrained optimization problem to achieve our design objective. The bio-inspired particle swarm method (PSO) is applied effectively to solve the optimization problem.

This study lays the groundwork for extending the methodology for more complex design problems in photonics. Examples include the design of two-dimensional periodic structures and optical resonators. Other areas to explore are plasmonics\cite{plasmonics} metals and in two-dimensional media\cite{tlow}.

\end{document}